\documentclass[12pt,reqno]{amsart}

\usepackage{amssymb,amsmath,graphicx,amsfonts,euscript}
\usepackage{color}

\newcommand{\be}{\begin{equation}}
\newcommand{\ee}{\end{equation}}

\newcommand{\ba}{\begin{array}{l}}
\newcommand{\ea}{\end{array}}

\newcommand{\na}{\nabla}

\numberwithin{equation}{section}

\begin{document}
\title[Inviscid Models]
{An incompressible 2D didactic model with singularity and explicit solutions of the
2D Boussinesq equations}

\author[D. Chae, P. Constantin and J. Wu]{Dongho Chae$^{1}$, Peter Constantin$^{2}$ and  Jiahong Wu$^{3}$}

\address{$^1$ Department of Mathematics,
College of Natural Science,
Chung-Ang University,
Seoul 156-756, South Korea}

\email{dchae@cau.ac.kr}

\address{$^2$ Program in Applied and Computational Mathematics,
Princeton University, Fine Hall, Washington Road, Princeton, NJ 08544-1000, USA.}

\email{const@math.princeton.edu}

\address{$^3$ Department of Mathematics,
Oklahoma State University,
401 Mathematical Sciences,
Stillwater, OK 74078, USA; and
Department of Mathematics,
Chung-Ang University,
Seoul 156-756, Republic of Korea}

\email{jiahong@math.okstate.edu}

\vskip .2in
\begin{abstract}
We give an example of a well posed, finite energy, 2D incompressible active scalar equation with the same scaling as  the
surface quasi-geostrophic equation and prove that it can produce finite time singularities. In spite of its simplicity, this seems to be the first such example. Further, we construct explicit solutions of the 2D Boussinesq equations whose gradients
grow exponentially in time for all time. In addition, we introduce  a variant of the 2D Boussinesq  equations which is perhaps a more faithful companion of the 3D axisymmetric Euler equations than the usual 2D Boussinesq equations.
\end{abstract}

\maketitle

\section{Introduction}

The purpose of this paper is threefold: first, to provide an example of an incompressible 2D active scalar, similar to the inviscid surface quasi-geostrophic (SQG) equation, that possesses a family of solutions which develop
finite-time singularities; second, to construct a class of explicit solutions to the inviscid
2D Boussinesq equations that grow exponentially in time; and third, to propose for study a
modified 2D Boussinesq system, that appears to provide a closer comparison to the
3D axisymmetric Euler equations than the standard 2D Boussinesq equations.
All the calculations are elementary, and the results serve as didactic examples.
\vskip .1in
The inviscid SQG equation is
\begin{equation} \label{SQG}
\left\{
\begin{array}{l}
\displaystyle \partial_t \theta + u\cdot\nabla \theta =0, \quad x \in \mathbb{R}^2, \, t>0,\\
\displaystyle  u =\nabla^\perp \psi \equiv (-\partial_{x_2}, \partial_{x_1}) \psi, \quad
  \Lambda \psi = \theta,
\end{array}
\right.
\end{equation}
where $\theta$ and $\psi$ are scalar functions of $x$ and $t$, $u$ denotes the 2D velocity field
and $\Lambda=(-\Delta)^{\frac12}$ denotes the Zygmund operator, which can be defined through the Fourier
transform
$$
\widehat{\Lambda f} (\xi)  = |\xi| \widehat{f}(\xi).
$$
The inviscid SQG equation is useful in modeling atmospheric phenomena
such as frontogenesis, the formation
of strong fronts between masses of hot and cold air (\cite{CMT,Gil,HPGS,Pe}). In addition,
the inviscid SQG equation is a significant example of a 2D
active scalar and some of its distinctive features have made it an important testbed for turbulence theories (\cite{Blu}). Mathematically, the inviscid SQG equation is difficult to analyze, and the issue of whether its solutions can develop singularities in a finite time remains open. Other active scalars with the same scaling $u\sim\theta$ include the porous medium (or Muskat) equation where
$u = (0, \theta) + \na p$. The blow up problem from smooth initial data is open there as well (see \cite{Cor} and references therein).
Here we give a first example of a 2D incompressible active
scalar equation with a velocity field having the same level of regularity and scaling as the active scalar, $u\sim\theta$,
\begin{equation} \label{SQGmodel}
\left\{
\begin{array}{l}
\displaystyle  \partial_t \theta + u\cdot\nabla \theta =0, \quad x \in \mathbb{R}^2, \, t>0,\\
\displaystyle  u =\nabla^\perp \psi, \quad
  -\partial_{x_2} \psi = \theta.
\end{array}
\right.
\end{equation}
The only difference between (\ref{SQG}) and (\ref{SQGmodel}) is the equation relating $\psi$ and $\theta$.
It is not very difficult to see that (\ref{SQGmodel}) is locally well-posed in a sufficiently regular
functional setting. What is more striking about this model is that (\ref{SQGmodel}) admits a family of
solutions which develop finite-time singularities even though they are initially smooth and have
finite energy.

\vskip .1in
Very recently Luo and Hou performed careful numerical simulations of the 3D axisymmetric incompressible Euler equations which suggested the appearance of a finite-time singularity \cite{LuoHou}. We briefly describe the setup and their
main result. The spatial domain is the cylinder
$$
\left\{(x_1,x_2, z): \,\, r\equiv \sqrt{x_1^2 + x_2^2} \le 1, \,\, 0\le z\le L  \right\}
$$
with periodic boundary conditions in the $z$-direction and no-penetration condition on the solid boundary $r=1$. The angular components of the velocity, vorticity and stream functions are
odd with respect to $z$ across $z=0$. The velocity on the unit circle $z=0$, $r=1$ vanishes, and thus all points on this circle are stagnation points. According to \cite{LuoHou}, the vorticity of a numerical solution at the
stagnation points  blows up in finite time. This finite-time singularity does not appear to be well understood theoretically. Motivated by the numerical simulations of Luo and Hou, Kiselev and Sverak recently proved the double exponential growth (in time) of the vorticity gradient of 2D Euler solutions in the unit disk \cite{KiSv}. The odd symmetry, the stagnation point and the boundary appear to be important in their work. Other pursuits for lower bounds for the vorticity gradient of the 2D Euler equation can be found in \cite{Denisov,Zlatos}. Our goal here is to provide explicit solutions to the 2D incompressible
Boussinesq equations which exhibit exponential growth in time.  The 2D Boussinesq equations are
\begin{equation} \label{Boussinesq}
\left\{
\begin{array}{l}
\displaystyle  \partial_t \omega + u\cdot\nabla \omega = \partial_{x_1} \theta,\quad x \in \mathbb{R}^2, \, t>0,\\
\displaystyle  \partial_t \theta + u\cdot\nabla \theta =0,\\
\displaystyle  u =\nabla^\perp \psi, \quad
  \Delta\psi = \omega.
\end{array}
\right.
\end{equation}
As pointed out in \cite{MB}, (\ref{Boussinesq}) is closely related to the 3D axisymmetric Euler equations.
We construct two families of solutions of (\ref{Boussinesq}). The first family is given in a fixed domain
with a stagnation point at the corner of the domain. The solutions in the family have an odd symmetry with respect to $x_2=0$, i.e. they are odd as functions of $x_2$, in the direction of gravity,
and  have $\nabla \theta$ growing exponentially in time. The second family of solutions
of (\ref{Boussinesq}) consists of smooth global solutions in a domain with a moving boundary.
Both $\nabla \omega$ and $\nabla \theta$ grow exponentially in time. It may be possible to further
exploit and extend these constructions to obtain solutions with gradients with double exponential growth.

\vskip .1in
We also propose for study the following modified 2D Boussinesq system
\begin{equation} \label{Boussinesq2}
\left\{
\begin{array}{l}
\displaystyle  \partial_t \omega + u\cdot\nabla \omega = -\partial_{x_2} (\theta^2),\quad x \in \mathbb{R}^2, \, t>0,\\
\displaystyle  \partial_t \theta + u\cdot\nabla \theta =0,\\
\displaystyle  u =\nabla^\perp \psi, \quad
  \Delta\psi = \omega.
\end{array}
\right.
\end{equation}
(\ref{Boussinesq2}) differs from (\ref{Boussinesq}) on the right-hand side of the vorticity equation, and has a different parity symmetry. Note that the gravity now points in the $x_1$ direction, which might be confusing but it is done to mimic the set-up of Luo and Hou, but the main difference with usual Boussinesq is that but both $\theta$ and $\omega$ are  allowed to be odd in $x_2$, i.e. in the direction perpendicular to gravity, not parallel to it.
As explained in Section \ref{mod}, (\ref{Boussinesq2}) appears to be a more exact match with the 3D axisymmetric
Euler equations as set up in the numerical simulations of Luo and Hou \cite{LuoHou}. We construct a class of
solutions of (\ref{Boussinesq2}) whose gradients exhibit exponential growth.

\vskip .1in
The rest of the paper is divided into three sections. The second section details the active scalar model similar to the SQG equation and describes its solutions which develop finite-time singularities. The
third section constructs two families of solutions to the 2D Boussinesq equations while the last section
presents a variant of the 2D Boussinesq equations and some of its explicit solutions with exponential gradient growth.

\vskip .3in
\section{A 2D model with finite-time singularity}
\label{SQGparallel}

This section presents an active scalar equation transported by an divergence-free
velocity field that admits solutions with finite-time singularities.

\vskip .1in
Consider the initial-value problem (\ref{SQGmodel}), namely
\begin{equation} \label{SQGIVP}
\left\{
\begin{array}{l}
\displaystyle  \partial_t \theta + u\cdot\nabla \theta =0, \quad x \in \mathbb{R}^2, \, t>0,\\
\displaystyle  u =\nabla^\perp \psi, \quad
  -\partial_{x_2} \psi = \theta,\\
\displaystyle  \psi(x,0)=\psi_0(x).
\end{array}
\right.
\end{equation}
It is not difficult to see that (\ref{SQGIVP}) is locally well-posed if $\psi_0\in H^s(\mathbb{R}^2)$ with $s>3$ and
$\theta_0=-\partial_{x_2} \psi_0\in L^1(\mathbb{R}^2)$. Now, we
assume that the initial stream function $\psi_0$, the initial active scalar $\theta_0$ and the
initial-velocity $u_0$ obey the following symmetry, for any $(x_1,x_2) \in \mathbb{R}^2$,
$$
\psi_0(x_1, -x_2) = -\psi_0(x_1, x_2), \qquad \theta_0(x_1, -x_2) = \theta_0(x_1, x_2),
$$
$$
u_{01}(x_1, -x_2) = u_{01}(x_1, x_2), \qquad u_{02}(x_1, -x_2) = -u_{02}(x_1, x_2).
$$
It is easy to check that the corresponding solution of (\ref{SQGIVP}) preserves this symmetry.
In particular, $u_2(x_1,0,t) =0$ and the equation for $\theta$ on the $x_1$-axis becomes
$$
\partial_t \theta(x_1,0,t) + u_1(x_1, 0,t) \partial_{x_1} \theta(x_1, 0,t) =0.
$$
or
\begin{equation}\label{Burgers}
\partial_t \theta(x_1,0,t) + \theta(x_1, 0,t) \partial_{x_1} \theta(x_1, 0,t) =0,
\end{equation}
which is the inviscid Burgers equation. It is well-known that Burgers' equation  can develop
discontinuities. In fact, if the initial data satisfies
$$
\partial_{x_1} \theta_0(x_1,0) <0 \quad\mbox{for some $x_1\in \mathbb{R}$},
$$
then the corresponding solution $\theta (x_1,0,t)$ of (\ref{Burgers}) becomes singular in
a finite time, namely $\partial_{x_1} \theta$ becomes infinity in a finite time.
One could also restrict to a periodic domain and obtain a finite time blow up. For example,
$$
\psi_0(x_1,x_2) =-\cos(x_1)\,\sin(x_2)
$$
will generate a solution that becomes singular in a finite time.

\vskip .3in
\section{Explicit solutions of the 2D Boussinesq equations}
\label{Bou}

In this section we construct two families of solutions to the 2D Boussinesq equations given by
(\ref{Boussinesq}). The first family of solutions is defined in a wedge and the gradient of $\theta$
grows exponentially in time. The second family of solutions is close to the first family, but it
is defined in a smooth domain with a moving boundary, and the gradients of both $\omega$ and
$\theta$ grow exponentially in time.

\vskip .1in
\subsection{Explicit solutions in a wedge}

The spatial domain is bounded by two half-lines in the positive half-plane:
$$
x_2 = 2x_1, \quad x_2 =-2x_1, \quad x_1 \ge 0.
$$
The stream function, the velocity field and the vorticity are given by
\begin{align*}
\psi(x) =& \left\{
\begin{array}{ll}
\frac{x_2^2}{2} - \,x_1 x_2, \quad & \,x_2 \ge 0, \\
-\frac{x_2^2}{2} - \,x_1 x_2, \quad & \,x_2 < 0;
\end{array}
\right. \\
u_1(x) =& -\partial_{x_2} \psi =\left\{
\begin{array}{ll}
-x_2 + x_1, \quad & \,x_2 \ge 0, \\
x_2 + x_1, \quad & \,x_2 < 0;
\end{array}
\right.
\\
u_2(x) =& \partial_{x_1} \psi = -x_2;
\\
\omega(x) =& \Delta \psi = \left\{
\begin{array}{ll}
1, \quad & \,x_2 \ge 0, \\
-1, \quad & \,x_2 < 0;
\end{array}
\right.
\end{align*}
It is clear that $\psi=0$ on the boundary of the domain and, consequently,
$u$ satisfies the no-penetration boundary condition. Given any initial data $\theta_0(x)$ depending only
on $x_2$ and odd with respect to $x_1$-axis, namely
$$
\theta_0(x) = \theta_0(x_2), \qquad \theta_0(-x_2) = -\theta_0(x_2), \qquad x_2\in \mathbb{R},
$$
the corresponding solution $\theta(x,t)$ preserves these properties and satisfies
\begin{equation} \label{rho}
\left\{
\begin{array}{l}
\displaystyle \partial_t \theta(x_2,t) - x_2 \, \partial_{x_2} \theta(x_2,t) =0, \\
\displaystyle \theta(x,0) =\theta_0(x_2).
\end{array}
\right.
\end{equation}
Integrating on characteristics, $\theta$ is given by
$$
\theta (x_2, t) = \theta_0(e^tx_2).
$$
with arbitrary odd $\theta_0$.  It is easy to see that the gradient of any solution to (\ref{rho}) grows exponentially in time. In fact,
$\partial_{x_2} \theta$ satisfies
\begin{align*}
\partial_t (\partial_{x_2} \theta) - \partial_{x_2} \theta - x_2 \partial_{x_2} (\partial_{x_2}\theta) =0,
\end{align*}
which is solved by
$$
\partial_{x_2} \theta = e^t\, (\theta_0')(e^tx_2).
$$
A special example is given by
\begin{equation*} \label{rhoexample}
\theta(x,t) = \sin(x_2\, e^t), \qquad \partial_{x_2} \theta(x,t) = e^t \cos(x_2\, e^t).
\end{equation*}
Many other explicit solutions can be obtained by taking different initial data.

\vskip .2in
\subsection{Explicit smooth solutions in a moving domain}

The spatial domain in the previous subsection has a sharp corner at the origin
and the vorticity has jumps when crossing the $x_1$-axis.  In this subsection we construct a smooth solution.

\vskip .1in
For $\sigma(x_2,t)$ to be explicitly determined later, we seek solutions with
stream function, velocity field and vorticity given by
\begin{align*}
\psi(x) =& \left\{
\begin{array}{ll}
\frac12 \sigma(x_2,t)\,x_2^2 - \,x_1 x_2, \quad & \,x_2 \ge 0, \\ \\
-\frac12 \sigma(x_2,t)\,x_2^2 - \,x_1 x_2, \quad & \,x_2 < 0;
\end{array}
\right.
\end{align*}
\begin{align*}
u_1(x) =& -\partial_{x_2} \psi =\left\{
\begin{array}{ll}
-\sigma(x_2,t)\,x_2 -\frac12 x_2^2\,\partial_{x_2} \sigma + x_1, \quad & \,x_2 \ge 0, \\ \\
 \sigma(x_2,t)\,x_2 + \frac12 x_2^2\,\partial_{x_2} \sigma + x_1, \quad & \,x_2 < 0;
\end{array}
\right.
\\
u_2(x) =& \partial_{x_1} \psi = -x_2;
\\
\omega(x) =& \Delta \psi = \left\{
\begin{array}{ll}
\frac12\, \partial_{x_2}^2(\sigma(x_2,t)\,x_2^2), \quad & \,x_2 \ge 0, \\ \\
-\frac12\, \partial_{x_2}^2(\sigma(x_2,t)\,x_2^2), \quad & \,x_2 < 0.
\end{array}
\right.
\end{align*}
Consequently, the boundary of the domain is given by $\psi=0$, namely
$$
2 x_1 = \sigma(x_2,t)\,x_2, \qquad 2 x_1 = -\sigma(x_2,t)\,x_2.
$$
Since the vorticity depends only on $x_2$,  the vorticity equation is reduced to
$$
\partial_t \omega - x_2 \, \partial_{x_2} \omega = 0, \qquad \omega(x_2,0) =\omega_0(x_2),
$$
which is solved by
$$
\omega(x,t) = \omega_0(e^tx_2).
$$
Therefore, $\sigma(x_2,t)$ must satisfy
$$
\partial_{x_2}^2(\sigma(x_2,t)\,x_2^2) =2 \omega_0(e^tx_2), \;{\mbox{for}}\; x_2>0\;{\mbox{and}}\;
-\partial_{x_2}^2(\sigma(x_2,t)\,x_2^2) =2 \omega_0(e^tx_2), \; {\mbox{for}}\; x_2<0.
$$
We also seek a solution $\theta$ depending on $x_2$ and $t$ only, namely $\theta=\theta(x_2,t)$, then
$$
\partial_t \theta - x_2 \, \partial_{x_2} \theta= 0, \qquad \theta(x_2,0) =\theta_0(x_2),
$$
which is solved by
$$
\theta(x,t) = \theta_0(e^tx_2).
$$
Some special examples of the solutions include
$$
\psi(x,t)= \frac16 x_2^3 \,e^t - x_1 x_2,\quad \omega = x_2 \,e^t, \quad \theta = x_2 \,e^t
$$
in the domain bounded by $x_1=\frac16 x^2_2\, e^t$.

\vskip .3in
\section{A variant of the 2D Boussinesq equations}
\label{mod}

In this section we suggest a new 2D Boussinesq system for study. This new system modifies
the standard 2D Boussinesq equations and appears to be a more direct parallel to the
3D axisymmetric Euler equations. A class of explicit solutions are constructed to exemplify the behavior of it solutions.

\vskip .1in
The modified 2D Boussinesq equations are given in (\ref{Boussinesq2}). We explain why this modified
system is identical to the 3D axisymmetric Euler equations, away from the axis of symmetry.
We first recall the 3D axisymmetric Euler equations in cylindrical coordinates (see, e.g., \cite{Chae,MB}). Let $\mathbf{e}_r$,
$\mathbf{e}_\theta$ and $\mathbf{e}_z$ be the orthonormal unit vectors defining the cylindrical coordinate
system,
$$
\mathbf{e}_r = \left(\frac{x_1}{r}, \frac{x_2}{r}, 0\right), \quad \mathbf{e}_\theta
=\left(-\frac{x_2}{r}, \frac{x_1}{r}, 0\right), \quad \mathbf{e}_z =(0, 0, 1),
$$
where $r=\sqrt{x_1^2 + x_2^2}$. A vector field $v$ is axisymmetric if
$$
v= v^r(r,z,t)\, \mathbf{e}_r + v^\theta(r,z,t)\, \mathbf{e}_\theta + v^z(r,z,t)\, \mathbf{e}_z.
$$
The 3D axisymmetric Euler equations can be written as
\begin{equation} \label{Euler}
\left\{
\begin{array}{l}
\displaystyle \frac{\widetilde{D}}{Dt} (r u^\theta) =0, \\ \\
\displaystyle \frac{\widetilde{D}}{Dt} \left(\frac{\omega^\theta}{r} \right) =-\frac1{r^4} \partial_z(ru^\theta)^2,\\ \\
\left(\partial^2_r + \frac1r \partial_r + \partial_z^2 -\frac1 {r^2}\right) \psi^\theta = \omega^\theta,
\end{array}
\right.
\end{equation}
where $u^\theta$, $\omega^\theta$ and $\psi^\theta$ are the angular components of the velocity,
vorticity and stream function, respectively, and
$$
\frac{\widetilde{D}}{Dt} =\partial_t + u^r \partial_r + u^z \partial_z.
$$
(\ref{Euler}) is a complete system since $u^r$ and $u^z$ can be recovered from $\psi^\theta$ by
$$
u^r = -\partial_z \psi^\theta, \qquad  u^z = \frac1r\, \partial_r (r \psi^\theta).
$$
Therefore the divergence-free condition
$$
\partial_r (r u^r) + \partial_z (r u^z) =0
$$
is automatically satisfied. We remark that the notation in the book of Majda  and Bertozzi \cite[p.62-66]{MB} is slightly
different from that in Luo and Hou \cite{LuoHou}. In particular, $\psi$ in \cite{MB} corresponds to
$r \psi^\theta$ in \cite{LuoHou} and $\omega^\theta$ in \cite{MB} corresponds to $-\omega^\theta$
in \cite{LuoHou}.

\vskip .1in
\cite{LuoHou} numerically solved (\ref{Euler}) in the cylinder
$$
D(1,L) = \left\{(r,z): \, 0\le r\le 1, \,\, 0\le z \le L \right\}
$$
with the initial data
$$
u_0^\theta(r,z) =100r\, e^{-30(1-r^2)^4}\,\sin\left(\frac{2\pi}{L}z\right),
\quad \omega^\theta_0(r,z) = \psi_0^\theta(r,z) =0
$$
subject to the periodic boundary condition in $z$,
$$
u^\theta(r,0,t) = u^\theta(r,L,t), \quad \omega^\theta(r,0,t) = \omega^\theta(r,L,t), \quad
\psi^\theta(r,0,t) = \psi^\theta(r,L,t)
$$
and the no-flow boundary condition on the solid boundary $r=1$,
$$
\psi^\theta(1,z,t) =0.
$$
Since the initial data $(u_0^\theta, \omega_0^\theta,\psi_0^\theta)$ is odd with respect to $z=0$,
the solution $(u^\theta, \omega^\theta,\psi^\theta)$ is also odd with respect to $z=0$. Consequently,
$$
u^r(1, z,t) =0, \quad u^z(r, 0, t) =0.
$$
In particular, at $z=0$ and $r=1$,
$$
u=(u^r, u^\theta, u^z) =0
$$
and all points on this circle are stagnation points.

\vskip .1in
The modified 2D Boussinesq system that we propose for study is given by (\ref{Boussinesq2}), namely
\begin{equation} \label{Boussinesq20}
\left\{
\begin{array}{l}
\displaystyle  \partial_t \omega + u\cdot\nabla \omega = -\partial_{x_2} (\rho^2),\quad x \in \mathbb{R}^2, \, t>0,\\ \\
\displaystyle  \partial_t \rho + u\cdot\nabla \rho =0,\\ \\
\displaystyle  u =\nabla^\perp \psi, \quad
  \Delta\psi = \omega.
\end{array}
\right.
\end{equation}
We changed the name of the scalar $\theta$ to $\rho$ in order to avoid confusing it with the notation for the angular variable.
We have the following correspondence between (\ref{Euler}) and (\ref{Boussinesq20}):
$$
 r \leftrightarrow x_1, \quad z\leftrightarrow x_2, \quad r u^\theta \leftrightarrow \rho,
 \quad \frac{\omega^\theta}{r} \leftrightarrow \omega, \quad (u^r, u^z) \leftrightarrow u,
 \quad \psi^\theta \leftrightarrow \psi.
$$
Thus, away from the symmetry axis, the behavior of the solutions to these two systems are
expected to be identical. In addition, both $\omega$ and $\rho$ can have the same parity symmetry with respect to $x_2=0$ (they are both odd). Note that the direction $x_2$ is not the direction of gravity, but it is perpendicular on it.

\vskip .1in
We seek a family of solutions to (\ref{Boussinesq20}) independent of $x_1$,
$$
\omega=\omega(x_2,t), \quad u_2=u_2(x_2,t), \quad \rho=\rho(x_2,t).
$$
The equations for $\omega$ and $\rho$ are reduced to
\begin{align*}
&\partial_t \omega + u_2 \partial_{x_2} \omega = -\partial_{x_2} (\rho^2), \\
&\partial_t \rho + u_2 \partial_{x_2} \rho =0.
\end{align*}
To make a comparison with the setup in \cite{LuoHou}, we further assume the odd symmetry with respect to
$x_2=0$. To simplify the calculation, we also assume that $u_2$ is time independent. From incompressibility and definitions, it is easy to see that $u_2$ is then a constant multiple of $x_2$, say
\[
u_2(x_2) = - x_2
\]
and consequently $\rho$ is given by
$$
\rho(x_2,t) = \rho_0(e^tx_2)
$$
and $\omega$ by
\[
\omega(x_2,t) = \omega_0(e^tx_2) - 2(e^t-1)\rho_0(e^tx_2)\rho_0'(e^tx_2)
\]
with $\rho_0, \omega_0$ arbitrary functions of one variable.
A special example is as follows:
\begin{align*}
\psi =& \left\{
\begin{array}{ll}
\left(1-\frac23e^t(e^t-1)x_2\right)\,\frac{x_2^2}{2} - \,x_1 x_2, \quad &\,x_2 \ge 0, \\ \\
-\left(1+\frac23e^t(e^t-1)x_2\right)\,\frac{x_2^2}{2} - \,x_1 x_2, \quad &\,x_2 < 0;
\end{array}
\right.
\end{align*}
\begin{align*}
u_1 =& -\partial_{x_2} \psi =\left\{
\begin{array}{ll}
-x_2 + e^t(e^t-1) x_2^2 + x_1, \quad & \,x_2 \ge 0, \\ \\
x_2 + e^t(e^t-1) x_2^2 + x_1, \quad & \,x_2 < 0;
\end{array}
\right.
\\
u_2 =&\partial_{x_1} \psi = -x_2
\end{align*}
and
\begin{align*}
\omega =& \Delta \psi = \left\{
\begin{array}{ll}
1-2x_2e^t(e^t-1), \quad & \,x_2 \ge 0, \\ \\
-1 - 2 x_2e^t(e^t-1), \quad & \,x_2 < 0;
\end{array}
\right.
\\
\rho =& x_2 e^t.
\end{align*}
One may also construct examples with fast oscillations such as
\begin{align*}
\omega =& \left\{
\begin{array}{ll}
1- (e^t-1)\sin(2x_2 e^t), \quad &\,x_2 \ge 0, \\ \\
-1 -(e^t-1)\sin(2x_2 e^t) , \quad & \,x_2 < 0;
\end{array}
\right.
\\
\rho =& \sin(2 x_2 e^t).
\end{align*}
Arbitrary growth (algebraic, exponential, etc) in $\rho_0$ in space engenders arbitrary growth in time in the solution (exponential, double exponential, etc).

\vskip .4in
\section*{Acknowledgements}
Chae was partially supported by NRF Grant no. 2006-0093854.
Constantin was partially supported by NSF grants DMS-1209394 and DMS-1265132.
Wu was partially supported by NSF grant DMS-1209153. Wu thanks Prof. C. Li, Prof. S. Preston and
Dr. A. Sarria for discussions. PC and JW thank the hospitality of the Department of
Mathematics of Chung-Ang University, where this work was performed.

\end{document}